\documentclass[11pt,twoside]{amsart}
\usepackage{amsmath,amssymb,amsfonts,enumerate,amsthm,graphicx,color}

\def\opn#1#2{\def#1{\operatorname{#2}}} 
\opn\chara{char} \opn\length{\ell} \opn\pd{pd} \opn\rk{rk}
\opn\injdim{inj\,dim} \opn\rank{rank} \opn\depth{depth}
\opn\grade{grade} \opn\height{height} \opn\embdim{emb\,dim}
\opn\codim{codim}

\opn\Tr{Tr} \opn\bigrank{big\,rank}
\opn\superheight{superheight}\opn\lcm{lcm}
\opn\trdeg{tr\,deg}%
\opn\reg{reg} \opn\lreg{lreg} \opn\skel{skel} \opn\set{set}
\opn\div{div} \opn\Div{Div} \opn\cl{cl} \opn\Cl{Cl}
\opn\Spec{Spec} \opn\Supp{Supp} \opn\supp{supp} \opn\Sing{Sing}
\opn\Ass{Ass}
\opn\Ann{Ann} \opn\Rad{Rad} \opn\Soc{Soc}
\opn\Ker{Ker} \opn\del{del} \opn\Im{Im} \opn\Hom{Hom} \opn\reg{reg}
\opn\Tor{Tor} \opn\Ext{Ext} \opn\End{End} \opn\Aut{Aut} \opn\id{id}

\opn\nat{nat}
\opn\pff{pf}
\opn\Pf{Pf} \opn\GL{GL} \opn\SL{SL} \opn\mod{mod} \opn\ord{ord}
\opn\aff{aff} \opn\con{conv} \opn\relint{relint} \opn\st{st}
\opn\lk{lk} \opn\cn{cn} \opn\core{core} \opn\vol{vol}
\opn\link{link} \opn\star{star} \opn\skel{skel} \opn\Reg{Reg}
\opn\gr{gr}

\def\pot#1#2{#1[\kern-0.28ex[#2]\kern-0.28ex]}

\opn\dirlim{\underrightarrow{\lim}}
\opn\inivlim{\underleftarrow{\lim}}
\def\Implies{\ifmmode\Longrightarrow \else
     \unskip${}\Longrightarrow{}$\ignorespaces\fi}
\def\implies{\ifmmode\Rightarrow \else
     \unskip${}\Rightarrow{}$\ignorespaces\fi}
\def\iff{\ifmmode\Longleftrightarrow \else
     \unskip${}\Longleftrightarrow{}$\ignorespaces\fi}

\let\:=\colon

\newtheorem{thm}{Theorem}[section]
\newtheorem{cor}[thm]{Corollary}
\newtheorem{lem}[thm]{Lemma}

\newtheorem{defn}[thm]{Definition}
\newtheorem{exam}[thm]{Example}
\newtheorem{rem}[thm]{\bf{Remark}}

\numberwithin{equation}{section}
\def\pn{\par\noindent}

\begin{document}


\title
{Bounds for the regularity of edge ideal of vertex decomposable and
shellable graphs}
\author{S. Moradi$^*$ and D. Kiani}

\thanks{{\scriptsize
\hskip -0.4 true cm MSC(2000): Primary: 13F55, 13D02; Secondary:
05C75
Received: 1 July 2009, Accepted: 23 October 2009\\
$*$Corresponding author
\newline\indent{\scriptsize $\copyright$ 2008 Iranian Mathematical
Society}}}

\maketitle

\begin{center}
Communicated by\;  Saeid Azam
\end{center}

\begin{abstract}
In this paper we give upper bounds for the regularity of edge ideal
of some classes of graphs in terms of invariants of graph. We
introduce two numbers $a'(G)$ and $n(G)$ depending on graph $G$ and
show that for a vertex decomposable graph $G$, $\reg(R/I(G))\leq
\min\{a'(G),n(G)\}$ and for a shellable graph $G$, $\reg(R/I(G))\leq
n(G)$. Moreover it is shown that for a graph $G$, where $G^c$ is a
$d$-tree, we have $\pd(R/I(G))=\max_{v\in V(G)} \{\deg_G(v)\}$.

\end{abstract}

\vskip 0.2 true cm


\pagestyle{myheadings}
\markboth{\rightline {\scriptsize  S. Moradi and D. Kiani}}
         {\leftline{\scriptsize Bounds for the regularity of edge ideal of vertex decomposable and shellable graphs}}

\bigskip
\bigskip


\vskip 0.4 true cm

\section{\bf Introduction}

\vskip 0.4 true cm Let $G$ be a simple graph with vertex set
$V(G)=\{x_1,\ldots,x_n\}$ and edge set $E(G)$. The edge ideal of $G$
in the polynomial ring $R=k[x_1,\ldots,x_n]$ is defined as
$I(G)=(x_ix_j:\{ x_i,x_j\}\in E(G))$. The edge ideal of a graph was
first considered by Villarreal \cite{V}. Finding connections between
algebraic properties of an edge ideal and invariants of graph is of
great interest. One question in this area is to explain the
regularity of an edge ideal by some information from graph. For some
classes of graphs for example chordal graphs  and shellable
bipartite graphs this question is answered, see \cite{HT1} and
\cite{VT}. For these graphs it is shown that the regularity of
$R/I(G)$ is equal to the maximum number of pairwise $3$-disjoint
edges in $G$, which is denoted by $a(G)$. Also in \cite[Lemma
2.2]{K}, it is shown that for any graph $G$, $\reg(R/I(G))\geq
a(G)$. In this paper we give upper bounds for $\reg(R/I(G))$ for
shellable and vertex decomposable graphs in terms of invariants of
graph. First we recall some definitions:

Let $G$ be a graph. An independent set of $G$ is a subset
$F\subseteq V(G)$ such that $e\nsubseteq F$, for any  $e\in E(G)$.
The \textbf{independence complex} of $G$ is the simplicial complex
$$\Delta_{G}=\{F\subseteq V(G): F\ \text{is an independent set of}\  G\}$$

For a simplicial complex  $\Delta$ on $X$ the \textbf{Alexander dual
simplicial complex} $\Delta^{\vee}$ to $\Delta$ is defined as
follows:

$$\Delta^{\vee}=\{F\subseteq X; X\setminus F\notin \Delta\}$$.

\begin{defn}
{\em A simplicial complex $\Delta$ is \textbf{shellable} if the
facets (maximal faces) of $\Delta$ can be ordered $F_1,\ldots,F_s$
such that for all $1\leqslant i<j\leqslant s$, there exists some
$v\in F_j\setminus F_i$ and some $l\in \{1,\ldots,j-1\}$ with
$F_j\setminus F_l=\{v\}$. We call $F_1,\ldots,F_s$ a shelling for
$\Delta$.}
\end{defn}
The above definition is referred to as {\it non-pure shellable} and
is due to Bj\"{o}rner and Wachs \cite{BW}. In this paper we will
drop the adjective ''non-pure". A graph $G$ is called shellable, if
the independence complex $\Delta_G$ is shellable.

\begin{defn}
{\em A monomial ideal $I=(f_1,\ldots,f_m)$ of the polynomial ring
$R=k[x_1,\ldots,x_n]$ has {\textbf {linear quotients}}, if there
exists an order $f_1<\cdots<f_m$ on the generators of $I$ such that
the colon ideal $(f_1,\ldots,f_{i-1}):f_i$ is generated by a subset
of variables for all $2\leq i\leq m$.}
\end{defn}

Also for any $1\leq i\leq m$, $\set_I(f_i)$ is defined as
$$\set_I(f_i)=\{x_k:\ x_k\in (f_1,\ldots, f_{i-1}) : f_i\}.$$
The following result relates squarefree monomial ideals with linear
quotients and shellable simplicial complexes:

\noindent{\textbf{Theorem A}} \cite[Theorem 1.4]{HD} The simplicial
complex $\Delta$ is shellable if and only if $I_{\Delta}^{\vee}$ has
linear quotients.

For a simplicial complex $\Delta$ and $F\in \Delta$, link of $F$ in
$\Delta$ is defined as $\lk_{\Delta}(F)=\{G\in \Delta: G\cap
F=\emptyset, G\cup F\in \Delta\}$ and the deletion of $F$ is the
simplicial complex $\del_{\Delta}(F)=\{G\in \Delta: G\cap
F=\emptyset\}$.
\begin{defn}
{\em Let $\Delta$ be a simplicial complex on the vertex set $V =
\{x_1,\ldots, x_n\}$. Then $\Delta$ is \textbf{vertex decomposable}
if either:

1) The only facet of $\Delta$ is $\{x_1,\ldots, x_n\}$, or
$\Delta=\emptyset$.

2) There exists a vertex $x\in V$ such that $\del_{\Delta}(x)$ and
$\lk_{\Delta}(x)$ are vertex decomposable, and such that every facet
of $\del_{\Delta}(x)$ is a facet of $\Delta$.}
\end{defn}
A graph $G$ is called vertex decomposable, if the independence
complex $\Delta_G$ is vertex decomposable.

The \textbf{Castelnuovo-Mumford regularity} (or simply regularity)
of an $R$-module $M$ is defined as: $$\reg(M) := \max\{j-i |\
\beta_{i,j}(M)\neq 0\},$$ and
$$\pd(M) := \max\{i |\ \beta_{i,j}(M)\neq 0 \ \text{for some}\ j\}.$$

For a monomial ideal $I=(x_{11}\cdots x_{1n_1},\ldots,x_{t1}\cdots
x_{tn_t})$ of the polynomial ring $R$, the \textbf{Alexande dual
ideal} of $I$ which is denoted by $I^{\vee}$ is defined as:
$$I^{\vee}=(x_{11},\ldots, x_{1n_1})\cap \cdots \cap (x_{t1},\ldots, x_{tn_t}).$$

 The following theorem was proved in \cite{T}.

\noindent{\textbf{Theorem B}}. Let $I$ be an square-free monomial
ideal. Then $\pd(I^{\vee})=\reg(R/I)$.

Two edges $\{x,y\}$ and $\{w,z\}$ of $G$ are called
\textbf{$3$-disjoint} if the induced subgraph of $G$ on
$\{x,y,w,z\}$ consists of exactly two disjoint edges or
equivalently, in the complement graph $G^c$, the induced graph on
$\{x, y,w, z\}$ is a four-cycle. A path of length $n$ is the graph
with $V(G)=\{x_1,\ldots,x_{n+1}\}$ and $E(G) = \{\{x_1, x_2\},
\{x_2, x_3\},\ldots, \{x_n, x_{n+1}\}\}$.

In this paper we find upper bounds for $\reg(R/I(G))$ in the case of
shellable and vertex decomposable graphs. In Theorem \ref{n}, we
show that for a shellable graph $G$, $\reg(R/I(G))\leq n(G)$ and in
Corollary \ref{corollary} it is shown that for a vertex decomposable
graph $G$, $\reg(R/I(G))\leq \min \{a'(G),n(G)\}$. In Theorem
\ref{1}, it is shown that if $G^c$ has no triangle, then
$\reg(R/I(G))\leq 2$ and finally Theorem \ref{pro} shows that for a
graph $G$ where $G^c$ is a $d$-tree, the projective dimention of
$R/I(G)$ is equal to $\max_{v\in V(G)}\{\deg_G(v)\}$.

\vskip 0.4 true cm

\section{\bf {\bf \em{\bf Main results}}}

\vskip 0.4 true cm
For a graph $G$, let $a'(G)$ be the maximum
number of vertex disjoint paths of length at most two in $G$ such
that paths of lengths one are pairwise $3$-disjoint in $G$. Also
from $\alpha'(G)$ we mean the matching number of $G$.

\vskip 0.4 true cm
\begin{thm}\label{sh}
Let $G$ be a vertex decomposable graph. Then $\reg(R/I(G))\leq
a'(G)$.
\end{thm}
\pn{\bf Proof.} By Theorem B, we have
$\reg(R/I(G))=\pd(I(G)^{\vee})$. So it is enough to show that
$\pd(I(G)^{\vee})\leq a'(G)$. By induction on $|V(G)|$ we prove the
assertion. For $|V(G)|=2$ there is nothing to prove. Let $|V(G)|>2$.
From the definition of vertex decomposable, there exists a vertex
$x\in V(G)$ such that $\del_{\Delta}(x)$ and $\lk_{\Delta}(x)$ are
vertex decomposable. Let $H_1=G\setminus \{x\}$ and $H_2=G\setminus
(\{x\}\cup N_G(x))$. It is easy to see that
$\del_{\Delta}(x)=\Delta_{H_1}$ and $\lk_{\Delta}(x)=\Delta_{H_2}$.
Thus $H_1$ and $H_2$ are vertex decomposable and each facet of
$\Delta_{H_1}$ is a facet of $\Delta_G$. Since a minimal vertex
cover of a graph is the complement of a facet of the independence
complex, for any minimal vertex cover $C$ of $H_1$, $C\cup\{x\}$ is
a minimal vertex cover of $G$. Also observe that for each minimal
vertex cover $C$ of $G$ containing $x$, $C\setminus \{x\}$ is a
minimal vertex cover of $H_1$. Therefore all the minimal vertex
covers of $G$ containing $x$ are $C_1\cup\{x\},\ldots,C_n\cup\{x\}$,
where $C_1,\ldots,C_n$ are the minimal vertex covers of $H_1$. Let
$N_G(x)=\{y_1,\ldots,y_t\}$ and let $C$ be a minimal vertex cover of
$G$ such that $x\notin C$. Then $\{y_1,\ldots,y_t\}\subseteq C$ and
$C\setminus \{y_1,\ldots,y_t\}$ is a minimal vertex cover of $H_2$.
Also for a minimal vertex cover $C$ of $H_2$,
$C\cup\{y_1,\ldots,y_t\}$ is a minimal vertex cover of $G$. Thus the
minimal vertex covers of $G$, which do not contain $x$ are
$C'_1\cup\{y_1,\ldots,y_t\},\ldots,C'_m\cup\{y_1,\ldots,y_t\}$,
where $C'_1,\ldots,C'_m$ are the minimal vertex covers of $H_2$.
Therefore $I(G)^{\vee}=xI(H_1)^{\vee}+y_1\cdots y_tI(H_2)^{\vee}$.
We show that $xI(H_1)^{\vee}\cap y_1\cdots
y_tI(H_2)^{\vee}=xy_1\cdots y_tI(H_2)^{\vee}$. Let $x^C\in
I(H_2)^{\vee}$ be a minimal generator. Then $C\cup\{y_1,\ldots,
y_t\}$ is a vertex cover of $H_1$ and hence $x^C\in I(H_1)^{\vee}$.
Thus $xy_1\cdots y_tI(H_2)^{\vee}\subseteq xI(H_1)^{\vee}\cap
y_1\cdots y_tI(H_2)^{\vee}$. Now let $x^C\in xI(H_1)^{\vee}\cap
y_1\cdots y_tI(H_2)^{\vee}$. Then $x,y_1,\ldots, y_t\in C$ and
$C\setminus \{x\}$ is a vertex cover of $H_1$ and
$C\setminus\{x,y_1,\ldots, y_t\}$ is a vertex cover of $H_2$. Thus
$x^C=xy_1\cdots y_tx^{C\setminus\{x,y_1,\ldots, y_t\}}\in xy_1\cdots
y_tI(H_2)^{\vee}$. Thus we have the following short exact sequence:

$$0\rightarrow xy_1\cdots
y_tI(H_2)^{\vee}\rightarrow xI(H_1)^{\vee}\oplus y_1\cdots
y_tI(H_2)^{\vee}\rightarrow I(G)^{\vee}\rightarrow 0.$$

Therefore $\pd(I(G)^{\vee})\leq
\max\{\pd(I(H_2)^{\vee})+1,\pd(I(H_1)^{\vee})\}$. By induction
hypothesis we have $\pd(I(H_1)^{\vee})\leq a'(H_1)$ and
$\pd(I(H_2)^{\vee})\leq a'(H_2)$. We consider two cases:

Case $1$. Let $\deg_G(x)\geq 2$, then $y_1,x,y_2$ is a path of
length two and $y_1,x,y_2\notin V(H_2)$. Thus $a'(H_2)+1\leq a'(G)$.
Since $a'(H_1)\leq a'(G)$, we have
$\pd(I(G)^{\vee})\leq\max\{a'(H_2)+1,a'(H_1)\}\leq a'(G)$.

Case $2$. Let $\deg_G(x)=1$ and $N_G(x)=\{y\}$ for some $y$. No
minimal vertex cover of $H_1$ contains $y$, since if a minimal
vertex cover of $H_1$ say $C$ contains $y$, then $C\cup\{x\}$ is a
non-minimal vertex cover of $G$, which is a contradiction as
discussed above. This means that each minimal vertex cover of $H_1$
contains $N_{H_1}(y)$. Thus $P_{N_{H_1}(y)}\subseteq\cap_{i=1}^n
P_{C_i}=I(H_1)$, where $P_{C_i}=(z: z\in C_i)$ and
$P_{N_{H_1}(y)}=(z: z\in N_{H_1}(y))$. Then $N_{H_1}(y)=\emptyset$,
since all the minimal generators of $I(H_1)$ are of degree two.
Therefore $x,y$ is a path which is $3$-disjoint from the paths of
length one in $H_2$ and disjoint from all paths in $H_2$. Thus
$a'(H_2)+1\leq a'(G)$. Since $a'(H_1)\leq a'(G)$, the assertion
follows from the inequality
$\pd(I(G)^{\vee})\leq\max\{a'(H_2)+1,a'(H_1)\}$.

\vskip 0.4 true cm
 H\`{a} and Van Tuyl in \cite{HT1} proved that
for any graph $G$, $\reg(R/I(G))\leq \alpha'(G)$, where $\alpha'(G)$
is the matching number, the largest number of pairwise disjoint
edges in $G$. It is easy to see that $a'(G)\leq \alpha'(G)$. The
following example shows that $a'(G)$ is a smaller upper bound for
vertex decomposable graphs.

\vskip 0.4 true cm

\begin{exam}
{\em  Let $G$ be a graph which is obtained from adding a vertex $x$
to the cycle $C_{2n+1}$ and joining it to one vertex of $C_{2n+1}$.
Let $y\in V(C_{2n+1})$ be a vertex that $xy\in E(G)$. Observe that
$H_1=G\setminus \{y\}$ and $H_2=G\setminus (\{y\}\cup N_G(y))$ are
path graphs and hence they are vertex decomposable. Also any facet
of $\Delta_{H_1}$ is a facet of $\Delta_{G}$. Therefore $G$ is
vertex decomposable. One can see that $\alpha'(G)=n+1$ and
$a'(G)=n$.}
\end{exam}

\vskip 0.4 true cm The following theorem was proved in \cite{cone}.

\begin{thm}\cite[Lemma 1.5]{cone}\label{herzog}
 Suppose that $I=(u_1,\ldots,u_m)$ is a monomial ideal with linear quotients with the ordering  $u_1<\cdots<u_m$ such that
$\deg(u_1)\leq \deg(u_2)\leq \cdots \leq \deg(u_m)$. Then the
iterated mapping cone $F$, derived from the sequence $u_1,\ldots,
u_m$, is a minimal graded free resolution of $I$, and for all $i>0$,
the symbols $$f(\sigma;u)\ \ \text{with}\ \ u\in G(I),\ \sigma
\subseteq \set_I(u),\ |\sigma|=i$$ form a homogeneous basis of the
$R$-module $F_i$. Moreover, $\deg f(\sigma;u)=|\sigma|+\deg(u)$.
\end{thm}
\vskip 0.4 true cm
 In the following theorem we show that for a shellable graph there
exists a vertex $x\in V(G)$ such that $\reg(R/I(G))$ is bounded by
$\reg(R/I(G\setminus(\{x\}\cup N_G(x)))+1$. For a subset $F\subseteq
V(G)$, the monomial $\prod_{x\in F}x$ is denoted by $x^F$.
\begin{thm}\label{shr}
Let $G$ be a shellable graph. There exists a vertex $x\in V(G)$ such
that if $H=G\setminus(\{x\}\cup N_G(x))$, then $$\reg(R/I(G))\leq
\reg(R/I(H))+1.$$
\end{thm}
\vskip 0.4 true cm \pn{\bf Proof.} By Theorem B, we have
$\reg(R/I(G))=\pd(I(G)^{\vee})$. Let $J=I(G)^{\vee}$. From Theorem
A, there exists an order of linear quotients $u_1<\cdots<u_t$ on the
minimal generators of $J$. From \cite[Lemma 2.1]{SZ}, one can assume
that $\deg(u_1)\leq \cdots \leq \deg(u_t)$. Thus by Theorem
\ref{herzog}, we have $\beta_i(J)=\sum_{j=1}^t{|\set_J(u_j)|\choose
i}$. Therefore $\pd(J)=\max\{|\set_J(u_i)|:\ 1\leq i\leq t\}$. For
any $i$, $1\leq i\leq t$, we have $u_i=x^{C_i}$, where $C_i\subseteq
V(G)$ is a minimal vertex cover of $G$. Let
$\pd(J)=|\set_J(x^{C_l})|$ for some $1\leq l\leq t$ and
$\set_J(x^{C_l})=(x_1,\ldots,x_r)$. Set $x=x_r$ and
$H=G\setminus(\{x\}\cup N_G(x))$ and $K=I(H)^{\vee}$. The set of
minimal vertex covers of $H$ is $\{C_i\setminus
N_G(x):N_G(x)\subseteq C_i \}$. Let $1\leq i_1<\cdots<i_k\leq t$ be
all integers such that $N_G(x)\subseteq C_{i_j}$ for $1\leq j\leq
k$. Then $K=(x^{C_{i_j}\setminus N_G(x)}: 1\leq j\leq k)$. Also the
ordering $x^{C_{i_1}\setminus N_G(x)}<\cdots<x^{C_{i_k}\setminus
N_G(x)}$ is an order of linear quotients for $K$ and it is degree
increasing. Since $x\in \set_J(x^{C_l})$, we have $x\notin C_l$.
Thus $N_G(x)\subseteq C_l$. Therefore $l=i_{l'}$ for some $1\leq
l'\leq k$. From the definition of linear quotients we see that for
any $1\leq i\leq r-1$, there exists $\lambda_i<l$ such that
$C_{\lambda_i}\setminus C_l=\{x_i\}$. It is easy to see that
$x\notin C_{\lambda_i}$ $(1\leq i\leq r-1)$. This means that
$N_G(x)\subseteq C_{\lambda_i}$ and consequently
$(x^{C_{\lambda_i}\setminus N_G(x)}:x^{C_l\setminus N_G(x)})=(x_i)$
for any $i$, $1\leq i\leq r-1$. Therefore $\set_K(x^{C_l\setminus
N_G(x)})=\{x_1,\ldots,x_{r-1}\}$. Thus $\reg(R/I(H))=\pd(K)\geq
|\set_K(x^{C_l\setminus N_G(x)})|=r-1=\reg(R/I(G))-1$.

\vskip 0.4 true cm Let $G$ be a graph and $x\in V(G)$. By a whisker
we mean adding a new vertex $y$ to $G$ and connecting $y$ to $x$.
This new graph is denoted by $G\cup W(x)$. We denote by $G\cup W(G)$
the graph obtained from $G$ by adding whiskers to all vertices of
$G$. In the following theorem the set of all induced subgraphs of
$G$ is denoted by $\mathcal{S}(G)$. \vskip 0.4 true cm
\begin{thm}\label{n}
Let $G$ be a shellable graph and $$n(G)=\max\{|V(H)|: H\in
\mathcal{S}(G), H\cup W(H)\in \mathcal{S}(G)\}.$$ Then
$\reg(R/I(G))\leq n(G)$.
\end{thm}
\vskip 0.4 true cm \pn{\bf Proof.} By Theorem B, it is enough to
show that $\pd(I(G)^{\vee})\leq n$. With the same notations as in
Theorem \ref{shr}, let $x^{C_1}<\cdots<x^{C_t}$ be an order of
linear quotients for $I(G)^{\vee}$ and
$\pd(I(G)^{\vee})=|\set_{I(G)^{\vee}}(x^{C_l})|=r$ for some $1\leq
l\leq t$. Let $\set_{I(G)^{\vee}}(x^{C_l})=(x_1,\ldots,x_r)$ and
$x^{C_{i_1}},\ldots,x^{C_{i_r}}<x^{C_l}$ be the monomials for which
$(x^{C_{i_j}}:x^{C_l})=(x_j)$ for any $1\leq j\leq r$. For any
$1\leq j\leq r$ we have $x_j\notin C_l$ and $x_j\in C_{i_j}$.
Therefore $N_G(x_j)\nsubseteq C_{i_j}$, since $C_{i_j}$ is a minimal
vertex cover of $G$. Also for any $1\leq j,k\leq r$, where $k\neq
j$, we have $x_k\notin C_{i_j}$, since $C_{i_j}\setminus C_l=
\{x_j\}$. For any $1\leq j\leq r$ let $y_j\in N_G(x_j)\setminus
C_{i_j}$. Thus for any $1\leq j,k\leq r$, where $k\neq j$ we have
$x_ky_j\notin E(G)$. Otherwise for the minimal vertex cover
$C_{i_j}$ we have $x_k\in C_{i_j}$ or $y_j\in C_{i_j}$, a
contradiction. Let $H$ be the induced subgraph of $G$ on the vertex
set $\{y_1,\ldots,y_r\}$. We have $x_jy_j\in E(G)$ and $x_jy_k\notin
E(G)$ for any $1\leq j,k\leq r$, where $k\neq j$. This means that
$H\cup W(H)\in \mathcal{S}(G)$. Therefore
$\pd(I(G)^{\vee})=r=|V(H)|\leq n$.

\vskip 0.4 true cm
\begin{exam}\label{masel}
{\em  Consider the graph $G$ with vertex set $\{x_1,\ldots,x_4\}$
and edge set $\{x_1x_2,x_1x_3,x_2x_3,x_1x_4\}$ Then $G$ is shellable
with $\reg(R/I(G))=1.$ We have $\alpha'(G)=2$ and $n(G)=1$. This
shows that $n(G)$ is an smaller upper bound for shellable graphs.}
\end{exam}
\vskip 0.4 true cm
\begin{exam}
{\em Let $G$ be a graph which is obtained from adding a vertex $x$
to the cycle $C_{2n+1}$ and joining it to two adjacent vertices of
$C_{2n+1}$. Then by \cite[Proposition 4.3]{DE}, $G$ is vertex
decomposable and hence it is shellable.  One can see that
$\alpha'(G)=n+1$. Obsereve that $n(G)\leq
\lfloor\frac{|V(G)|}{2}\rfloor=n+1$. We show that $n(G)<n+1$. By
contradiction assume that $n(G)=n+1$. Let $H$ be an induced subgraph
of $G$ such that $n(G)=|V(H)|=n+1$. Then $|H\cup W(H)|=2n+2$. Hence
$H\cup W(H)=G$. Thus $G$ has $n+1$ vertices of degree one, a
contradiction. Therefore $n(G)<\alpha'(G)$.}
\end{exam}
\vskip 0.4 true cm
\begin{rem}
{\em There are graphs for which $a'(G)<n(G)$. The path graph of
length three is such an example for which $a'(G)=1$ and $n(G)=2$.
Also there are graphs for which $n(G)<a'(G)$. Consider the complete
graph $K_n$ for $n\geq 6$. We have $n(G)=1$ and $a'(G)\geq 2$.}
\end{rem}
\vskip 0.4 true cm
\begin{cor}\label{corollary}
Let $G$ be a vertex decomposable graph. Then $\reg(R/I(G))\leq
\min\{a'(G),n(G)\}$.
\end{cor}
\vskip 0.4 true cm \pn{\bf Proof.} This follows from Theorems
\ref{sh}, \ref{n} and the fact that every vertex decomposable graph
is shellable, which is proved in \cite[Theorem 11.3]{BW}. \vskip 0.4
true cm
\begin{thm}\label{1}
Let $G$ be a graph such that $G^c$ has no triangle, then
$\reg(R/I(G))\leq 2$. In addition if $G^c$ is not chordal, then
$\reg(R/I(G))=2$.
\end{thm}

\vskip 0.4 true cm \pn{\bf Proof.} From Hochster's formula we have
$$\beta_{i,j}(R/I(G))=\sum_{S\subseteq V;|S|=j} \dim
\widetilde{H}_{j-i-1}(\Delta(G^c_S),K),$$ where $G_S$ denotes the
induced subgraph of $G$ on the vertex set $S$. Since $G^c$ has no
cycle of length $3$, any clique in $G^c$ is of cardinality at most
$2$. Thus $\widetilde{H}_i (\Delta(G^c_S),K)=0$ for any $i>1$ and
any $S$. Therefore $\widetilde{H}_{j-i-1}(\Delta(G^c_S),K)=0$ for
any $j-i>2$. Thus for any $i$ and $j$ such that
$\beta_{i,j}(R/I(G))\neq 0$, one has $j-i\leq 2$ and the result
holds. If $G^c$ is not chordal, then by \cite[Theorem 1]{F}, $I(G)$
does not have a linear resolution and hence $\reg(R/I(G))\neq 1$.
Thus $\reg(R/I(G))=2$.

\vskip 0.4 true cm
\begin{defn}
{\em A $d$-tree is a chordal graph defined inductively as follows:

$(i)$ $K_{d+1}$ is a $d$-tree.

$(ii)$ If $H$ is a $d$-tree, then so is $G=H\cup_{K_d} K_{d+1}$.}
\end{defn}

\vskip 0.4 true cm Edge ideals with $2$-linear resolution are
characterized in \cite[Theorem 1]{F} and it is shown that $I(G)$ has
linear resolution precisely when $G^c$ is a chordal graph. Eliahou
and Villarreal in \cite{EV} conjectured that $\pd(R/I(G))$, where
$I(G)$ has $2$-linear resolution, is equal to the maximum degree of
vertices of $G$. In the following theorem we show that for a graph
$G$ such that $G^c$ is a $d$-tree, we have $\pd(R/I(G))=\max_{v\in
V(G)} \{\deg_G(v)\}$. This statement is not true for an arbitrary
ideal with $2$-linear resolution. Consider the cycle $C_4$. Clearly
$C_4^c$ is chordal, and hence $I(C_4)$ has $2$-linear resolution but
$\pd(R/I(C_4))=3$, while $\max_{v\in V(C_4)} \{\deg_{C_4}(v)\}=2$.

To prove Theorem \ref{pro} we need the following easy lemma. \vskip
0.4 true cm
\begin{lem}\label{degree}
Let $G$ be a $d$-tree. Then $\deg_G(v)\geq d$ for any $v\in V(G)$.
\end{lem}

\vskip 0.4 true cm \pn{\bf Proof.} We proceed inductively in terms
of the definition of a $d$-tree. If $G=K_{d+1}$, then the assertion
is clear. Let $G=H\cup_{K_d}K_{d+1}$, where $H$ is a $d$-tree. Then
by induction hypothesis $\deg_H(v)\geq d$ for any $v\in V(H)$. Let
$V(G)=V(H)\cup\{x\}$, where $\{x\}=V(K_{d+1})\setminus V(H)$. Then
$\deg_G(x)=d$ and for any $v\in V(H)$, we have $\deg_G(v)\geq
\deg_H(v)\geq d$.

\vskip 0.4 true cm
\begin{thm}\label{pro}
Let $G$ be a graph such that $G^c$ is a $d$-tree. Then
$\pd(R/I(G))=\max_{v\in V(G)} \{\deg_G(v)\}$.
\end{thm}

\vskip 0.4 true cm \pn{\bf Proof.} We prove by induction on $|V(G)|$
that $I(G)$ has linear quotients and $\pd(R/I(G))=\max_{v\in V(G)}
\{\deg_G(v)\}$. For $|V(G)|=2$ the result is clear. Let $|V(G)|>2$
and $G'=G^c$. Here we have $G'=H\cup_{K_d} K_{d+1}$, where $H$ is a
$d$-tree. Let $V(G')\setminus V(H)=\{x\}$ and $V(H)\cap
V(K_{d+1})=\{x_1,\ldots,x_d\}$ and $V(H)\setminus
V(K_d)=\{y_1,\ldots,y_k\}$. Since $H$ is a $d$-tree, by induction
hypothesis $I(H^c)$ has linear quotients and
$\pd(R/I(H^c))=\max_{v\in V(H^c)} \{\deg_{H^c}(v)\}$. We have
$I(G)=(xy_1,\ldots,xy_k)+I(H^c)$. Let $u_1<\cdots<u_l$ be an order
of linear quotients for the minimal generators of $I(H^c)$. We claim
that the ordering $xy_1<\cdots<xy_k<u_1<\cdots<u_l$ is an order of
linear quotients for $I(G)$. Consider two monomials $xy_i$ and $u_j$
for some $1\leq i\leq k$ and $1\leq j\leq l$ and let $u_j=zw$ for
some $z,w\in V(H)$. Since $\{z,w\}$ is not an edge of $H$, then at
least one of $z$ and $w$ is not in $V(K_d)$. Without loss of
generality assume that $w\notin V(K_d)$. Then $w=y_{j'}$ for some
$1\leq j'\leq k$. We have $x|(xy_i:u_j)$ and $xy_{j'}<u_j$ and
$(xy_{j'}:u_j)=(x)$. For $xy_i<xy_j$, we have $(xy_i:xy_j)=(y_i)$
and for $u_i<u_j$, since $u_1<\cdots<u_l$ is an order of linear
quotients, the result holds. Now by Theorem \ref{herzog}, we have
$\pd(I(G))=\max\{|\set_{I(G)}(zw)|:\ \{z,w\}\in E(G)\}$ and
$\pd(I(H^c))=\max\{|\set_{I(H^c)}(zw)|:\ \{z,w\}\in E(H^c)\}$. For
any $1\leq i\leq k$, we have
$\set_{I(G)}(xy_i)=\{y_1,\ldots,y_{i-1}\}$. For any $1\leq j\leq l$,
we know that $u_j=y_{j'}z_j$ for some $1\leq j'\leq k$ and some
$z_j\in V(H)$. Thus $\set_{I(G)}(u_j)=\{x\}\cup \set_{I(H^c)}(u_j)$.
Therefore $\pd(I(G))=\max\{\pd(I(H^c))+1,k-1\}$ and hence
$\pd(R/I(G))=\pd(I(G))+1=\max\{\pd(R/I(H^c))+1,k\}$. Since
$\pd(R/I(H^c))=\max_{v\in V(H^c)} \{\deg_{H^c}(v)\}$, thus
$$\pd(R/I(G))=\max_{v\in V(H^c)}\{\deg_{H^c}(v)+1,k\}.$$ For any $i$,
$1\leq i\leq k$, we have $\deg_{G}(y_i)=\deg_{H^c}(y_i)+1$, because
$x$ is adjacent to $y_i$ in $G$. We claim that for any $1\leq i\leq
d$, $\deg_{G}(x_i)<k$. Let $1\leq i\leq d$ be an integer. Since $H$
is a $d$-tree, by Lemma \ref{degree} we have $\deg_H(x_i)\geq d$. So
there exists $y_j$ for some $1\leq j\leq k$ such that $x_iy_j\in
E(H)$. Therefore $\deg_{H^c}(x_i)<k$. Thus $\deg_{H^c}(x_i)+1\leq k$
for any $i$, $1\leq i\leq d$. Since $\deg_{H^c}(x_i)=\deg_{G}(x_i)$,
then $\deg_{G}(x_i)<k$ for any $i$, $1\leq i\leq d$ . Since
$\deg_{G}(x)=k$, thus $\max_{v\in
V(H^c)}\{\deg_{H^c}(v)+1,k\}=\max_{v\in V(G)} \{\deg_{G}(v)\}$ and
the proof is complete.


\vskip 0.4 true cm

\begin{center}{\textbf{Acknowledgments}}
\end{center}
The research of the second author was in part supported by a grant
from IPM with number (No. 88050116). The authors express their deep
gratitude to Professor Siamak Yassemi for his helpful suggestions
for the improvement of this work. They also thank the referees for
their useful comments. \\ \\


\bigskip
\bigskip

{\footnotesize \pn{\bf S. Moradi}\; \\ {Department of Pure
Mathematics,
 Faculty of Mathematics and Computer Science,
 Amirkabir University of Technology,
424, Hafez Ave., Tehran 15914, Iran.}\\ {\tt
Email:s\_moradi@aut.ac.ir}
{\footnotesize \pn{\bf D. Kiani}\; \\
{  Department of Pure Mathematics,
 Faculty of Mathematics and Computer Science,
 Amirkabir University of Technology,
424, Hafez Ave., Tehran 15914, Iran \&  School of Mathematics,
Institute for Research in Fundamental Sciences (IPM), P.O. Box
19395-5746, Tehran, Iran. }\\
{\tt Email:dkiani@aut.ac.ir}\\

\begin{thebibliography}{20}

\bibitem{BW} A. Bj\"{o}rner, M. L. Wachs, Shellable nonpure complexes and
posets I, {\em Trans. Amer. Math. Soc.} {\bf 348} (1996) no. 4,
1299–-1327.

\bibitem{DE} A. Dochtermann, A. Engstr\"{o}m,  Algebraic properties of edge ideals via
combinatorial topology, arXiv:0810.4120v1.

\bibitem{EV} S. Eliahou, R. Villarreal,   The second Betti number of an
edge ideal, {\em Aportaciones Matematicas, Serie Comunicaciones}
(1999) {\bf 25}:115–-119.

\bibitem{F} R. Fr\"{o}berg, On Stanley-Reisner rings, {\em Topics in Algebra} {\bf 26} (1990).

\bibitem{HT} H. T. H\`{a}, A. Van Tuyl,  Splittable ideals and the resolutions of monomial ideals, {\em J.
Algebra}
{\bf 309} (2007) 405--425.

\bibitem{HT1} H. T. H\`{a}, A. Van Tuyl,  Monomial ideals, edge ideals of hypergraphs,
and their graded Betti numbers, {\em J. Algebraic Combin.} {\bf 27}
(2008) no. 2, 215–-245.

\bibitem{HD} J.~Herzog, T.~Hibi and X.~Zheng,  Diracs theorem on chordal graphs and Alexander duality, {\em European J. Combin.} {\bf 25} (2004)
949--960.

\bibitem{cone} J. Herzog, Y. Takayama, Resolutions by mapping cones. {\em  Homology, Homotopy  Appl.} {\bf 4} (2002) no.
2, part 2, 277-–294.

\bibitem{K} M. Katzmann,  Characteristic-independence of Betti numbers of graph ideals, {\em J. Combin. Theory
Ser. A } {\bf 113} (2006)  no. 3, 435–-454.

\bibitem{SZ} A. Soleyman Jahan, X. Zheng,  Pretty clean monomial
ideals and linear quotients, arXiv:0707.291.4v1.
\bibitem{T} N. Terai,  Alexander duality theorem and Stanley-Reisner
rings, {\em S\"{u}rikaisekikenky\"{u}sho K\"{o}ky\"{u}ruko} (1999)
no. 1078, 174--184, Free resolutions of coordinate rings of
projective varieties and related topics (Kyoto 1998).

\bibitem{VT} A. Van Tuyl,  Sequentially Cohen-Macaulay bipartite graphs: vertex decomposability
and regularity, arXiv:0906.0273v1.


\bibitem{VV}
A.~Van Tuyl and R.~Villarreal,  Shellable graphs and sequentially
Cohen-Macaulay bipartite graphs, {\em J. Combin. Theory Ser. A} {\bf
115} (2008) 799--814.


\bibitem{V} R.H. Villarreal,  Cohen-Macaulay graphs, {\em  Manuscripta Math} {\bf 66} (1990)
no. 3, 277-–293.

\end{thebibliography}
\end{document}